\documentclass[12pt]{article}
\title{Some remarks on the algebra of bounded Dirichlet series}
\date{ April 15,  2008}
\begin{document}
\author{ Brian Maurizi, Herv\'e   Queff\'elec}
\maketitle

\newtheorem{definition}{Definition}[section]
\newtheorem{theo}{Theorem}[section]
\newtheorem{prop}{Proposition}[section]
\newtheorem{cor}{Corollary}[section]
\newtheorem{lem}{Lemma}[section]
\newtheorem{rem}{Remark}
\newtheorem{ex}{Example}
\def\etc{{\it etc}}
\def\im{\mathrm{Im}\,}
\def\ie{i.e.\ }
\def\tq{\ ;\ } 
\def\ffi{\varphi}
\def\eps{\varepsilon}
\def\noi{\noindent}
\newcommand {\dis}{\displaystyle}	
\newcommand {\N}{{\bf{N}}} 
\newcommand {\Z}{{\bf{Z}}}		
\newcommand {\R}{{\bf{R}}} 
\newcommand {\Q}{{\bf{Q}}} 
\newcommand {\C}{{\bf{C}}} 
\newcommand {\D}{\bf{D}} 
\newcommand {\T}{{\bf{T}}}
\newcommand {\F}{\bf{F}}  
\newcommand {\K}{\bf{K}}  
\newcommand {\B}{\bf{Bad}}


\section{ Introduction} The aim of this paper is to contribute to the study of the algebra of bounded Dirichlet series. But we must first recall several definitions, notations, and facts.
The analytic theory of Dirichlet series is similar to that of power series, but with important differences: whereas a power series has {\it one} radius of convergence, a Dirichlet series

$$  f(s)= \sum_1^\infty a_n n^{-s}, s=\sigma +it\eqno(1.1)$$
has {\it several} abscissas of convergence, we list four of them:
$$\begin{array}{lllll}\sigma_c&=&inf\{a:(1.1)\textrm {converges for}\  \Re s>a\}\\&=&\textrm {abscissa of simple convergence}\\ 
$$\sigma_u&=&inf\{a:(1.1)\textrm{converges uniformly for}\ \Re s>a\}\\&=&\textrm{abscissa of uniform convergence}\\
$$\sigma_a&=&inf\{a:(1.1)\textrm {converges absolutely for}\ \Re s>a\}\\&=&\textrm{abscissa of absolute convergence}\\
$$\sigma_b&=&inf\{a:(1.1) \textrm{has an analytic,  bounded extension for}\  \Re s>a\}\\&=&\textrm {abscissa of boundedness}\ .\end{array} $$

It is easy to see that $$\sigma_a -\sigma_c \leq1, \sigma_b\leq\sigma_u\leq\sigma_a$$
and the example $a_n=e^{in^{\alpha}}, 0\leq\alpha\leq 1$, for which $\sigma_c=1-\alpha, \sigma_a=1$(\cite {KAQ})
shows that the difference $\sigma_a -\sigma_c$ can take any value between $0$ and $1$. Now, in 1913, H.Bohr(\cite {BOHR 1}) proved the following theorem, which appears to be {\it basic} in the context of bounded Dirichlet series:

\begin{theo}: If the sum of a Dirichlet series $\sum_1^\infty a_n n^{-s}$, convergent for some value of $s$, has an analytic and bounded extension in a vertical half-plane, then the series converges uniformly in each smaller half-plane. In other words, one has $\sigma_b=\sigma_u$.
 \end{theo}
 H.Bohr then naturally asked three questions:
 
\smallskip
 
\noi{\bf Question 1 :} We have $\sigma_b=\sigma_u$; have we even $\sigma_b=\sigma_a$, i.e. have we $\sigma_a= \sigma_u$?

\smallskip
 
\noi{\bf Question 2 :} If not, what is the supremum $T$ of all possible differences $\sigma_a -\sigma_u $? (clearly, $0\leq T\leq1 $).

\smallskip
 
\noi{\bf Question 3 :} What about the absolute convergence of bounded {\it power} series?
 (\cite{BOHR 3}).
To tackle these questions ( now solved), Bohr was led to consider, through the Kronecker Approximation Theorem, Dirichlet series as Taylor series in infinitely many complex variables. But to formulate theorems more precisely, and in a more modern language, it will be convenient to introduce some symbols:
 
 $\C_\theta$ denotes the open half-plane $\{ s:\Re s>\theta\}$, where $\theta\in \R$, and $\T$ denotes the unit circle of the complex plane.
 
 ${\cal H}^\infty$ denotes the set of Dirichlet series $f(s)=\sum_1^\infty a_n n^{-s}$, analytic and bounded in $\C_0$, equipped with the sup-norm $$\Vert f\Vert_\infty =\sup_{s\in\C_0}\vert f(s)\vert$$
 which makes it a (unital and commutative) Banach algebra.
 
If $ n=p_1^{\alpha_1}...p_r^{\alpha_r}$ is a positive integer decomposed in prime factors  ( $p_1=2, p_2=3$, etc..., and $\alpha_j\geq 0$), $\Omega (n)=\alpha_1+...\alpha_r$ denotes the number of prime divisors of n counted with their multiplicity, and $ P^+(n)$ denotes the largest prime divisor of $n$.
If $d$ is a positive integer, ${\cal H}^{\infty}(d)$ denotes the Banach subspace (not a subalgebra) of ${\cal H}^\infty$ formed by those $f(s) =\sum_1^\infty a_n n^{-s} \in {\cal H}^\infty$ such that 
$$a_n\not = 0 \Rightarrow \Omega (n)\leq d.$$
For example, $${\cal H}^{\infty}(1) =\{\sum a_p p^{-s}\} ,\ {\cal H}^{\infty}(2)=\{\sum a_{pq} (pq)^{-s}\},$$
where p, q are primes.
$H^\infty$ will denote the (non-separable) Banach algebra of functions  which are analytic and bounded in the open unit disk $\D$, equipped with the norm of the supremum on $\D$. The map $ f\to f(2^{-s})$ is an isometry of $H^\infty$ onto a ( proper) subalgebra of ${\cal H}^\infty$, and in particular the latter algebra is non-separable.
Finally, $c_0$ denotes the Banach space of sequences tending to zero at infinity, with its natural norm, and $B=c_0\bigcap\D^\infty$ the open unit ball of that Banach space.
We also recall the 
\begin{theo}(Kronecker Approximation Theorem): Let $\lambda_1,...,\lambda_r$ be rationally independent real numbers (i.e. $\sum_1^r n_j\lambda_j=0\textrm\  {and}\ n_j\in \Z \Rightarrow
n_1=...n_r=0$).Then, the map
$$t\Rightarrow (e^{i\lambda_{1}t},...,e^{i\lambda_{r}t}):\R\to \T^r$$
has dense range.

\end{theo}
We might take $\lambda_j =\log p_j$, due to the uniqueness of the expansion of a positive integer in primes.
Now,a good deal of energy was invested to find the exact value of $T$, which turned out to be ${1\over 2}$; and this paper is accordingly  organized as follows:

 In Section 1 ( this introduction), we recall the "Bohr point of view", according to which the algebra ${\cal H}^\infty$ is viewed as a space of functions in an infinite number of complex variables.  In other words, we can view  ${\cal H}^\infty$ in two different ways:

\begin{enumerate}
 \item As an algebra of  functions on the half-plane $\C_0$
 \item As an algebra of functions on the the infinite-dimensional ball $B$ of $c_0$
\end{enumerate}

The first perspective leads to a new and simple proof of Bohr's Theorem 1.1, through a Theorem recently proved by the second-named author and al. on the behaviour of partial sums of bounded Dirichlet series. The second perspective leads (Theorem 1.4) , in a very simple way, to the failure of the Corona Theorem for the Banach algebra ${\cal H}^\infty$ on $\C_0$, a result which appears to be new.  Section 2 is mainly expository, it recalls the solutions to Questions 1,2,3 by Bohnenblust -Hille and others, as well as some recent refinements. Section 3 contains the two main results of the paper (Ths.3.1 and 3.2), either the statements or the proofs ( through a new deterministic device) being new.
 
 We will be more specific  in a moment on  the so-called  "Bohr's point of view"(\cite {BOHR 2}). Here is now  a typical example of the first perspective, which was proved in\cite{BCQ}:
 \newpage
\begin{theo}:Let $$f(s) = \sum_1^\infty a_n n^{-s}\in {\cal H}^\infty,\textrm{and}\      S_N(f)=\sum_1^N a_n n^{-s}.$$ 
Then, we have:
$\Vert S_N(f)\Vert_\infty\leq C\log N\Vert f\Vert_\infty$, where $N\geq 2$ and $C$ is a numerical constant.\end{theo}
The proof  uses classical tools of one -complex variable theory, like the Cauchy and Perron's formulas. And this Theorem    gives a new and  simple proof of Bohr's theorem 1.1, of which it is a quantitative version: Indeed, we have to show that, for each $\epsilon >0$, the series $\sum_1^\infty a_n n^{-s -\epsilon}$ converges uniformly in $\C_0$; now, setting $S_n= S_n(s)$, an Abel summation by parts gives 
$$\sum_1^N a_n n^{-s-\epsilon}=\sum_1^Nn^{-\epsilon}(S_n -S_{n-1})=\sum_1^{N-1}
S_n(n^{-\epsilon} -(n+1)^{-\epsilon})+N^{-\epsilon}S_N$$
and the series $\sum_1^\infty S_n(n^{-\epsilon} - (n+1)^{-\epsilon})$ is normally convergent in the half-plane $\C_0$ since its general term is dominated by ${\epsilon C \log n\over n^{\epsilon +1}}\Vert f\Vert_\infty$, whereas $\vert N^{-\epsilon}S_N(s)\vert \leq {C\log N\over N^{\epsilon}}\Vert f\Vert_\infty$, which proves Theorem 1.1.

The second perspective was discovered by Bohr(\cite {BOHR 2}), and can be summarized as follows: let $f(s)= \sum_1^\infty a_n n^{-s}$ be a "Dirichlet polynomial"( i.e. $a_n = 0$ for $n$ large, say $\geq N$), so that $$f\in {\cal H}^\infty \textrm \ {and}\   \Vert f\Vert_\infty= \sup_{t\in R}\vert f(it)\vert$$ (by the maximum modulus principle). Let $r=\pi(N)$, i.e. $p_r\leq N<p_{r+1}$. Then, each $n\leq N$ has a unique expansion $ n= p_1^{\alpha_1}...p_r^{\alpha_r}$, with $\alpha_j = \alpha_{j}(n)$. If $z= (z_1,...,z_r)\in \T^r$, or $z\in \D^r$, we set with Bohr:
$$ \Delta f(z)=\sum_1^N a_n z_1^{\alpha_{1}}...z_r^{\alpha_{r}}.\eqno(1.2)$$
Let $m$ be the Haar measure of $\T^r$. Bohr(\cite {BOHR 2}) observed that, as a consequence of the Kronecker Approximation Theorem and of the distinguished maximum principle, we have:
$$ \Vert f\Vert_\infty = \Vert \Delta f\Vert _\infty;  \Vert f\Vert_2=\Vert \Delta f\Vert_2.\eqno(1.3)$$
In the first equation, the sup norm of $\Delta f$ refers either to $\T^r$ or to the polydisk $\D^r$, while in the second equation $\Vert \Delta f\Vert_2$ refers to the Haar measure m, and $\Vert f\Vert_2$ to the Haar measure of the Bohr compactification of $\R$, i.e.
$$ \Vert f\Vert_2^{2}= \lim_{T\to \infty} {1\over 2T} \int _{-T}^T \vert f(it)\vert^2 dt= \sum_1^\infty \vert a_n\vert^2.\eqno(1.4)$$
It can be proved(\cite{HLS}) that this procedure can be extended to arbitrary functions f of ${\cal H}^\infty $ as follows : for $z=(z_1,z_2,...,z_r,...) \in B$, denote by $z^{(m)}$ the truncated sequence $(z_1,z_2,...,z_m,0,...,0,...)$ and define 
$\Delta f(z^{(m)})$ in the same way as we did for polynomials. It follows from the Schwarz lemma that , for $ z\in B$ and $l<m$, one has
$$\vert \Delta f(z^{(m)})- \Delta f(z^{(l)})\vert \leq 2\  \max_{l<j\leq m}\vert z_j\vert\Vert f\Vert_\infty,$$ so that $\Delta f (z^{(m)}$ tends to a limit $\Delta f(z)$ as $m$ tends to $\infty$, and we still have
$$ \Vert \Delta f\Vert_\infty= \Vert f\Vert_\infty,\eqno((1.5)$$
so that we have a homomorphism $\Delta:{\cal H}^\infty\to H^\infty(B)$ which is norm-preserving and in particular injective. As a corollary of that point of view, H.Bohr proved the following inequalities for $f(s) =\sum_1^\infty a_n n^{-s}\in {\cal H}^\infty$ ($p$ denoting a prime)
$$   (\sum_1^\infty  \vert a_n\vert^2)^{1/2}\leq \Vert f\Vert_\infty,\ \textrm{and}\ \sum\vert a_p\vert \leq \Vert f \Vert _\infty.\eqno(1.6)$$
As observed by H.Bohr, the first inequality in (1.6) easily implies that $T\leq 1/2$ in Question 2 ( this can also be proved without appealing to the Bohr point of view), which turned out to be the right value of $T$. Another example of the second perspective is given by the following theorem
\begin{theo} 1) The invertible elements of ${\cal H}^\infty$ are the functions which are bounded below.

2) There exist $f_1, f_2 \in {\cal H}^\infty$ such that $\vert f_1(s)\vert +\vert f_2(s)\vert \geq \delta>0$ for any $s\in \ \C_0$ and yet $f_1 g_1+f_2 g_2\not =1$ for any pair of functions $g_1, g_2 \in {\cal H}^\infty$.
\end {theo}
{\bf Proof :} 1) If $\vert f(s)\vert \geq \delta $, letting $\Re s$ tend to infinity, we get $\vert\ a_1\vert\geq \delta$, where $a_1$ is the constant term of $f$ , and then  an application of Neumann's lemma  shows that $1/f$ can be expanded as a Dirichlet series if $\Re s$ is large enough.  Hewitt and Williamson (\cite{HW}) elaborated on this point, proving furthermore that, if $f$ has an absolutely convergent Dirichlet series, so has ${1\over f}$. 

2) Suppose that $f_1 g_1 +f_2 g_2 = 1$ on $\C_0$. Set $F_1=\Delta f_1$ and define similarly $F_2,G_1,G_2$. As we saw before, the initial equation  implies that
 
$$  F_1(z) G_1(z) +F_2(z) G_2(z) =1 \quad  \forall z\in B.\eqno(1.7)$$
In particular, we must have, for some $\delta>0$
$$ \vert F_1(z)\vert +\vert F_2(z)\vert \geq \delta \  \forall z\in B.\eqno(1.8)$$ 
But it is easy to produce examples where the assumption $\vert f_1(s)\vert +\vert f_2(s)\vert \geq \delta>0$ holds, and not (1.8). Take for example $$f_1(s) = {1\over 2}+2^{-s}, \textrm{and}\  f_2(s)=3^{-s}.$$
Separating the cases $\Re s>2$ and $\Re s \leq 2$, we see that $\vert f_1(s)\vert +\vert f_2(s)\vert \geq 1/9$, but yet (1.8) fails, and even $F_1(z)={1\over2}+z_1 , \ F_2(z)= z_2$ have the common zero $z=(-{1\over 2}, 0,...0,...)\in B$! In other terms, the celebrated Corona Theorem of Carleson(\cite {CA}) for ${\cal H}^\infty$ considered as an algebra on $\C_0$ fails( Therefore , it also fails for absolutely convergent Dirichlet series). Here, the right point of view is that of infinitely many complex variables, and the true (open) question is to know if the necessary condition   (1.8) is sufficient to imply the existence of a Bezout identity $f_1g_1+f_2 g_2=1$.

\medskip

\noi{\bf Remark :} In\cite {NI}, the following criterion for the existence of  Bezout identities of length 2 (say) for a uniform algebra $A$ of bounded functions on a given set $X$  is used:
\begin{prop} The following are equivalent:\ 
1) Any element  $\varphi \in SpA$ is $2$-visible from X, i.e : for any pair $f_1,f_2$ of functions of A and any number $\epsilon >0$ , there exists some $x\in X$ such that:
$$\vert f_1(x) -\varphi (f_1)\vert <\epsilon;\ \vert f_2(x)-\varphi(f_2)\vert <\epsilon.$$
2) For any pair $f_1,f_2$ of functions in A such that $$\vert f_1(x)\vert +\vert f_2(x)\vert \geq \delta >0\ \forall x\in X,$$ there exist $g_1,g_2\in A$ such that $f_1g_1+f_2g_2=1$.
\end{prop}
As a confirmation of Theorem 1.2, we can exhibit quite a few elements of $Sp {\cal H}^\infty$ which are $2$-invisible from $\C_0$ under the form of the following (easy to prove) Proposition, in which $\chi$ denotes a completely multiplicative function( $\chi(mn)=\chi(m)\chi(n)\quad\forall m,n)$ defined on the set of positive integers and such that
$$\sum_1^\infty \vert a_n\ \chi(n)\vert\ < \infty \  \forall f(s)=\sum_1^\infty a_n n^{-s}\in {\cal H}^\infty.$$
\begin{prop} Let $\varphi(f) =\sum_1^\infty a_n\chi(n)\in SpA$. Then, the following are equivalent:

1) $\varphi$ is $2$-visible from $\C_0$.

2) $\vert \chi(n)\vert =n^{-c}$ for some $c\geq {1\over 2}$, and possibly $c=\infty$.

3) $\varphi$ is in the closure of $\C_0$ in $SpA$.
\end {prop}
This proposition shows that all completely multiplicative functions which fail to verify 2) will give rise to characters which are $2$-  invisible from $\C_0$.

As concerns question 3), H.Bohr( \cite {BOHR 3}) proved the following:
\begin{theo}
 Let $f(z)=\sum_0^\infty a_n z^n\in H^\infty$. Then:

$$\sum_0^\infty \vert a_n\vert ({1\over 3})^n\leq \Vert f\Vert_\infty.$$
Moreover, ${1\over 3}$ is optimal: if $\sum_0^\infty \vert a_n\vert r^n\leq \Vert f\Vert_\infty$ for each $f\in H^\infty$, we must have $r\leq {1\over 3}$.
\end {theo}
It is interesting to note that Theorem 1.4 was proved in the course of trying to answer Questions 1 and 2. For ${1\over3}<r<1$, we have
$$\sum_0^\infty \vert a_n\vert r^n \leq \ (\sum_0^\infty \vert a_n\vert^2)^{1\over 2}(\sum_0^\infty r^{2n})^{1\over 2} =(1-r^2)^{-1\over 2}\Vert f\Vert_2\leq (1-r^2)^{-1\over 2}\Vert f\Vert_\infty,$$
and one can ask about the best constant $C_r$ such that $\sum_0^\infty \vert a_n\vert \leq C_r \Vert f\Vert _\infty$ for each $f\in H^\infty$. Bombieri(\cite{BOM}) computed the exact value of $C_r$ for ${1\over 3}\leq r\leq {1\over \sqrt 2}$, which is $C_r= {1\over r}(3-\sqrt{1-r^2})$ and later Bombieri and Bourgain(\cite{BOMB}) proved that
$$ (1.9)\  C_r<(1-r^2)^{-1\over2}\  \forall r>{1\over \sqrt 2};\  C_r\sim (1-r^2)^{-1\over 2}$$
as $r\to 1$. We will see that the situation is quite different for the the space 
 ${\cal H}^\infty$ of bounded Dirichlet series. 
    
\section{The Hille-Bohnenblust Theorem and some refinements}

We recall that, with the notations of section 1,\  $T$ is the supremum of all possible differences $\sigma_a -\sigma_u=\sigma_{a}(f) -\sigma_{u}(f)$ as $f$ runs through ${\cal H}^\infty$, and that $0\leq T\leq {1\over 2}$. Similarly, we could define, for an integer $d$:
$$T_d=\sup\{\sigma_{a}(f)-\sigma_{u}(f)\}, f\in {\cal H}^\infty(d),$$
where ${\cal H}^\infty(d)$was defined in section 1. In a celebrated paper published in the Annals( \cite {BOHN}), Bohnenblust and Hille proved the following
\begin{theo} We have the following values for $T$ and $T_d$:

a) $T={1\over 2}$; b) $T_d={1\over 2}-{1\over {2d}}$. Moreover, $T$ and $T_d$ are attained.
\end {theo}
Note that in particular we have $T_1=0$, which is the content of inequality (1.6). The authors moreover proved that $\sigma_a -\sigma_u$ can take any value between $0$ and ${1\over 2}$ and Bohr later gave a very simple proof of that fact.
Now, the Bohnenblust-Hille Proof had two main ingredients:

1) An intensive use of Bohr's point of view that ${\cal H}^\infty$ is embedded in $H^\infty(B)$ as described before.

2) A clever use of Walsh matrices; we will return to those matrices.

This work was "prophetic" in some respects: it announced the Rudin -Shapiro sequence(see for example \cite{RU}), and the theory of p-Sidon sets in Harmonic analysis, through the use of Littlewood's inequality for multilinear forms; but it is fair to say that this work remained somehow qualitative; later, probabilistic proofs replaced the deterministic one of the two previous authors. What we shall do in  Section 3 will be  to revisit those two approaches, deterministic and probabilistic, under a more quantitative form, in order to get refined versions of Theorem 2.1. In this Section, we shall first  reformulate the definition of $T$ and $T_d$ in finite terms, and recall a probabilistic Lemma to be used in Section 3, as well as recent refinements of the Bohnenblust-Hille Theorem: 
\begin {theo} We have the following:
$$ T=inf\{\sigma\geq 0:\sum_1^N\vert a_n\vert\leq C_{\sigma}N^\sigma\Vert \sum_1^N a_n n^{-s}\Vert_\infty,\forall N,  a_1,...a_N\}.\eqno(2.1)$$
$$\textrm{ Same for } T_d, \textrm{except that}\   a_n=0 \ \textrm{whenever} \ \Omega (n)>d.\eqno(2.2)$$
\end {theo}
{\bf Proof:} Let $E$ be the set of $\sigma's$ in the Right hand Side  of (2.1), and $\sigma>T$.
If $f(s) =\sum_1^\infty a_n n^{-s}\in {\cal H}^\infty$, we know from Bohr's Theorem that $\sigma_{u}(f)\leq 0$, therefore $\sigma_{a}(f)\leq T<\sigma$, and $\sum_1^\infty\vert a_n\vert n^{-\sigma}<\infty$.
The closed graph theorem implies the existence of a finite  constant $C_\sigma$ such that
$$\sum_1^\infty \vert a_n\vert n^{-\sigma}\leq C_\sigma \Vert\sum_1^\infty a_n n^{-s}\Vert_\infty,\ \forall\  \sum_1^\infty a_n n^{-s}\in {\cal H}^\infty.$$
In particular, we have for given N:
$\sum_1^N\vert a_n\vert n^{-\sigma}\leq C_\sigma\Vert \sum_1^Na_n n^{-s}\Vert_\infty$, so that
$$\sum_1^N\vert a_n\vert\leq N^\sigma\sum_1^N\vert a_n\vert n^{-\sigma}\leq C_\sigma N^\sigma \Vert \sum_1^Na_n n^{-s}\Vert_\infty.$$
This implies that $\sigma\in E$, showing that $]T,\infty[\subset E$ and that 
$\inf E\leq T$. Conversely, suppose $\sigma\in E$. Let $$f(s) =\sum_1^\infty a_n n^{-s}\in {\cal H}^\infty,\ A_n=\sum_1^N\vert a_j\vert,\ \alpha>0.$$
An Abel's summation by parts shows that
$$\sum_1^N{\vert a_n\vert \over n^\alpha}=\sum_1^{N-1} A_n(n^{-\alpha}-(n+1)^{-\alpha})+A_N N^{-\alpha}\leq C_\alpha \sum_1^\infty {\log n\over n^{\alpha-\sigma+1}}<\infty$$
if $\alpha >\sigma$ (we used  Theorem 1.3). This shows that $\alpha\geq \sigma_{a}(f)\geq \sigma_{a}(f) -\sigma_{u}(f)$, and since $f$ is arbitrary:
$\alpha\geq T$. Finally, $T\leq\inf E$, and $T=\inf E$, as claimed. The proof of (2.2) is the same. \hfill$\diamondsuit$.

The form of (2.1) suggests the following: let ${\cal P_{N}}$ be the set of Dirichlet polynomials $P(s) =\sum_1^N a_n n^{-s}$,  $\Vert P\Vert_W=\sum_1^N\vert a_n\vert$ be the Wiener norm of $P$ and $S(\Lambda_N)$  (the Sidon constant of the set $\Lambda_N=\{\log1,...,\log N\}) $ be  the quantity
$$  S(\Lambda_N)=\sup_{P\in {\cal P_N}}{\Vert P\Vert _W\over \Vert P\Vert_\infty}\leq \sqrt N.\eqno(2.3)$$
What is the behaviour of $S(\Lambda_N)$ as $N\to \infty?$
In\cite {Q} , we proved the following:
\begin{theo} There exist numerical constants $a_0,b_0$ such that
$$  S(\Lambda_N)\geq a_0\sqrt N exp(-b_0 \lambda(N)),\eqno(2.4) $$
where $\lambda(x) =\sqrt{\log x\log\log x}$ for $x$ positive and large.We can take 
$b_0=\sqrt 2$.
\end {theo}
In view of (2.1), the previous theorem shows that $T\geq {1\over 2}$, that is to say it contains the Bohnenblust-Hille Theorem. The proof of (2.4) was based on the following probabilistic Lemma (see\cite{KA} , \cite{Q} or \cite {KOQ}), as well as on sharp estimates in the Rankin problem of Number Theory:
\begin{lem} let $P^(s) =\sum_1^N a_n n^{-s}\in {\cal P_N}$, such that:
$a_n\not =0\Rightarrow P^{+}(n)\leq y$, where $y$ is a fixed real number $\geq 2$. 
Let $$P_{\omega}(s)=\sum_1^N a_n\epsilon_{n}(\omega)n^{-s},$$
where $(\epsilon_n)$ is a sequence of i.i.d. Rademacher variables defined on some probability space $(\Omega,{\cal A},P)$,i.e. $P(\epsilon_{n}=1) =P(\epsilon_{n}=-1)= {1\over 2}$. Then, ${\cal E}$ denoting expectation, we have:
$$  {\cal E}\Vert P_{\omega}\Vert _\infty \leq C\Vert P\Vert_2\sqrt{{y\over \log y}}\sqrt {\log\log N}.\eqno(2.5)$$
\end{lem}
In\cite{KOQ}, we proved a converse of Theorem 2.3:
\begin{theo} There exist numerical constants $a_1,b_1>0$ such that
$$  S(\Lambda_N)\leq a_1\sqrt N exp(-b_1\lambda(N)).\eqno(2.6)$$
\end {theo}
Recently, R. de la Bret\`eche(\cite{LAB}) , still using lemma 2.1 and probabilistic techniques, improved on Theorems 2.3 and 2.4, showing that we can take $b_0={\sqrt 2\over 2},\  b_1={\sqrt 2\over 4}$. Independently, but still with random methods, K.Seip(\cite {SE}) obtained $b_0=1$. In\  \cite {BCQ} , we used a combination of Theorems 1.3 and 2.4 to prove that:
\begin{theo}
a) Let $f(s)=\sum_1^\infty a_n n^{-s}\in {\cal H}^\infty$.Then we have:
$$\sum_1^\infty \vert a_n\vert n^{-{1\over 2}}\leq C\Vert f\Vert_\infty.$$

b) If moreover $f\in {\cal H}^\infty (d)$ , we have:
$$\sum_1^\infty \vert a_n\vert n^{-\sigma_{d}}(\log n)^{{d-1\over 2}}\leq C\Vert f\Vert_\infty,$$
where $\sigma_d={1\over2} -{1\over 2d}$.
\end{theo}
Observe that Theorem 2.5 shows that there is no analogy of the extension by Bombieri-Bourgain of the Bohr Theorem 1.5: here, for $\sigma >{1\over 2}$, we have 
$$\sum_1^\infty \vert a_n\vert n^{-\sigma}\leq \sqrt {\sum_1^\infty\vert a_n\vert^2}\sqrt{\zeta(2\sigma)},\eqno(2.7)$$
and we might expect, as in Bombieri and Bourgain, that this is essentially optimal and that things explode as $\sigma\to {1\over 2}$, whereas a) and b) above show that this is not the case.
In the next Section, we will introduce a deterministic method, which can be viewed both as an extension of the Bohnenblust-Hille construction and as an extension of the Rudin-Shapiro sequence; this will allow us to reprove, without probabilities, Theorem 2.3 with $b_0=1$ , and to show that the exponent ${d-1\over 2}$ of Theorem 2.5 is optimal, a result which will be also given by lemma 2.1.

\section{An extension of the Bohnenblust-Hille and Rudin-Shapiro devices}
This Section contains the main results( Theorems 3.1 and 3.2 below) of the paper.
Let us first recall some facts; a {\it Walsh} matrix $A=(a_{ij})$ of  size $q$ is a square $q\times q$ matrix with unimodular coefficients ($\vert a_{ij}\vert = 1$) such that $A^\ast A= qI$, where $I$ is the identity matrix of size $q$. We can take $a_{ij}=\gamma_i(x_j)$, where $\{x_1,...,x_q\}$ is an abelian group $G$ of order $q$ and $\Gamma =\{\gamma_1,...,\gamma_q\}$ its dual; if $G$ is cyclic and $\omega$ is a primitive $q^{th}$ - root of unity, we can take $a_{ij}=\omega^{ij}$, and $A$ is called a {\it Schur} matrix; if $q=2^r$ and if $G$ is the "$2^r$-Gruppe of Klein", then $A$ has $\pm 1$-valued entries and is called a {\it Hadamard} matrix (such matrices also exist for other numbers than powers of $2$, e.g. for $q=20$). Hadamard matrices can be generated by blocks, according to the following inductive relation
$$ A_1=\pmatrix{1&1\cr1&-1\cr}, \quad A_{k+1}=\pmatrix{A_k&A_k\cr A_k&-A_k\cr}.$$

Observe that $A_1$ is the matrix appearing in the parallelogram law
$$\vert a+b\vert^2 +\vert a-b\vert^2= 2(\vert a\vert^2 +\vert b\vert^2), $$
which reads $\Vert A_1v\Vert^2= 2\Vert v\Vert^2$ if $v$ is the vector of $\C^2$ with coordinates $a$ and $b$. The proof of Bohnenblust and Hille makes use of Schur matrices; and later Shapiro(1951) and Rudin(1959) (see e.g.(\cite{RU}), independently, considered the following sequence $(P_n,Q_n)$ of pairs of polynomials:
$P_0=Q_0 =1$;\  then, $$P_{n+1}=P_n +z^{2^n}Q_n,\ Q_{n+1}=P_n - z^{2^n}Q_n.$$
For example, $$ P_1(z)=1+z,\  Q_1(z)=1-z;\ P_2(z)=1+z+z^2-z^3;\ Q_2(z) =1+z-z^2+z^3.$$
$P_n$ and $Q_n$ have $\pm 1$-valued coefficients. Now, the construction of that sequence can be described as an alternation of {\it shifts} and of {\it actions of $A_1$}:
$$(P_0,Q_0)\buildrel shift\over {\longrightarrow} (P_0,zQ_0)\buildrel A_1\over {\longrightarrow}  (P_1,Q_1)\longrightarrow\cdots$$
$$    (P_n,Q_n)\buildrel shift\over {\longrightarrow}    (P_n,z^{2^n}Q_n)\buildrel A_1\over {\longrightarrow}  (P_{n+1},Q_{n+1})\longrightarrow\cdots$$
the shift having for effect of making the spectra of $P_n$ and $z^{2^n} Q_n$ non-overlapping and the action of $A_1$ having for effect to keep moduli under control:
$$\vert z\vert =1\Longrightarrow \vert P_n\vert^2+\vert Q_n\vert^2=2(\vert P_{n-1}\vert^2+\vert Q_{n-1}\vert^2)=...=2^n(\vert P_0\vert^2+\vert Q_0\vert^2) =2^{n+1},$$
so that $\Vert P_n\Vert_\infty\leq 2^{n+1\over 2}$, whereas $\Vert P_n\Vert_W=2^n $.
If one tries to imitate that construction for Dirichlet polynomials in view of minorizing the Sidon constant $S(\Lambda_N)$, one is naturally led to the following:
$P_0= Q_0=1$; then, $$P_{n+1}(s) =P_{n}(s) +p_{n+1}^{-s} Q_n(s), Q_{n+1}(s)=P_{n}(s) -p_{n+1}^{-s}Q_n(s).$$
Set $N= p_1...p_n.$ We get in the same way:
$$P_n(s) =\sum_1^N a_k k^{-s},\vert a_k\vert =0\textrm \  {or}\  1\ ,\Vert P_n\Vert_W =2^n,\ \Vert P_n\Vert_\infty \leq 2^{{n+1\over 2}},$$
so $S(\Lambda_N)\leq 2^{{n-1\over 2}}$; this shows that $S(\Lambda_N)\to \infty$ with $N$, but this is all : as $n\to \infty$, $N\geq n!\geq n^{n} e^{-n}$, and the minoration $S(\Lambda_N)\geq 2^{{n-1\over 2}}$ gives the existence of uniformly ( in $\C_0$) and non-absolutely convergent Dirichlet series, but it does not even give the positivity of $T$ ( see (2.1)), the "degree" of $P_n$ being too large, compared with its "size" $2^n$. We therefore try to make the "degree" and "size" vary independently. To that effect , we {\it  pass to the Bohr point of view}, trying to produce algebraic polynomials with many variables , of comparatively low degree and with a big size, the two parameters  {\it degree} and {\it size} being more or less at our disposal. Fix a large integer $q\geq 1$, and define inductively sequences of $q-$ tuples of polynomials, by applying alternatively shifts (in independent variables ) and the action of a Walsh matrix $A=(a_{ij})$ of size $q$. Set $P_0^{(1)}=...=P_0^{(q)} =1$. Then proceed as follows:
$$( P_0^{(1)},...,P_0^{(q)})\buildrel shift\over {\longrightarrow}(z_1 P_0^{(1)},...,z_q P_0^{(q)})\buildrel A\over {\longrightarrow}(P_1^{(1)},...,P_1^{(q)})$$ 
$$\buildrel shift\over {\longrightarrow} (z_{q+1} P_1^{(1)},...,z_{2q}P_1^{(q)})
\buildrel A\over {\longrightarrow} (P_2^{(1)},...,P_2^{(q)}){\longrightarrow}\cdots$$
$$   
(P_d^{(1)},..., P_d^{(q)}) \buildrel shift\over {\longrightarrow} (z_{dq+1}P_d^{(1)},...,z_{(q+1)d} P_d^{(q)})\buildrel A\over {\longrightarrow} (P_{d+1}^{(1)},...,P_{d+1}^{(q)}).$$
For example:
$$P_1^{(i)}=\sum_{j=1}^q a_{ij_1}z_{j_1}; P_2^{(i)}=\sum_{j=1}^q z_{q+j_2}a_{ij_2} P_1^{(j_2)}=\sum_{j_2=1}^q z_{q+j_2} a_{ij_2}\sum_{j_1=1}^q  a_{j_{2}j_{1}} z_{j_1}.$$

The following flexible lemma summarizes the properties of the sequence of polynomials thus obtained:
\begin {lem}:
Fix $1\leq i\leq q$, and an integer $d\geq 1$. Then, the polynomial $P_{d}^{(i)}=P$ has the following properties:

a) The degree of P is $d$, and $P$ is homogeneous;

b) the number of variables in $P$ is $r=qd;$

c) $\Vert P\Vert_W = q^d;$

d) $\Vert P\Vert _\infty \leq q^{{d+1\over 2}};$

e) $P=\Delta Q$, where $$Q=\sum_1^N a_n n^{-s},\vert a_n\vert  =0\textrm\ {or}1,N\geq p_{r}^d \textrm\  {and}\ a_n\not =0\Longrightarrow n\textrm\ {squarefree \ and}\  \Omega (n)\leq d,$$ i.e. we have $Q\in {\cal H}^\infty (d).$
\end {lem}
{\bf Proof:}  a) At each inductive step, the degree is increased by one ( multiplication by $z_{qd+j}$ when one passes from $d$ to $d+1$).

b) At each inductive step, we add q new variables $z_{qd+1},z_{qd+2},...,z_{qd+q}$.

c) At each step, the size is multiplied by $q$, since we add independent variables, and the initial size is $1=q^0$. Observe ( cf.$P_2^{(i)} =\sum a_{ij_2} a_{j�{2}j_{1}}z_{j_{1}} z_{j_{2}+q}$) that $P$ has unimodular coefficients, so that $\Vert P\Vert_W =q^d =\textrm{size of}\  P$.

d) Suppose all the $z_j$ unimodular, as we may. Then, if $z=(z_1,...,z_q,...)$, the Walsh character of $A$ gives:
$$\sum_1^q\vert P_{d}^{(i)}(z)\vert^2 =q\sum _1^q\vert z_{(d-1)q +j}P_{(d-1)q}^{(j)}(z)\vert^2 =q\sum_1^q\vert P_{(d-1)}^{(j)}(z)\vert ^2$$ $$=...=q^d\sum_1^q\vert P_0^{(j)}(z)\vert ^2 =q^{d+1}.$$
e) A monomial of $P$ is of the form $u z_{j_1} z_{q+j_2}...z_{(d-1)q+j_d}$, where 
$\vert u\vert =1$ and $ 1\leq j_1,...,j_d\leq q$, so that $Q(s) = \sum_1^Na_n n^{-s}$, with $$\vert a_n\vert =0\textrm\ {or}\ 1,\textrm {and}\  a_n\not =0\Longrightarrow n=p_{i_1} p_{i_2}...p_{i_d}$$ with:

$$ 1\leq i_1\leq q; q+1\leq i_2\leq 2q;...q(d-1)+1\leq i_d\leq qd =r.\eqno(3.1)$$
In particular, the largest integer  in the spectrum of $Q$ is less than  $p_r^d$ .\hfill$\diamondsuit$

As a first corollary of the previous construction, we have a deterministic proof of the optimality of b) in Theorem 2.5, as well as a very simple proof of the Bohnenblust -Hille Theorem 2.1. We will denote by $\alpha$ the supremum of those exponents such that
$$ \sum_1^\infty \vert a_n\vert n^{-\sigma_{d}} (\log n)^\alpha \leq C_\alpha \Vert f\Vert _\infty,\ \forall f\in {\cal H}^\infty (d).\eqno(3.2)$$ 
From Theorem 2.5, we know that $\alpha\geq {d-1\over 2}$.
\begin{lem} The number $\alpha$ is the supremum of those exponents $\beta$ such that
$$\sum_1^N \vert a_n\vert\leq C_{\beta}{N^{\sigma_d}\over (\log N)^\beta}\Vert f\Vert_\infty,\ \forall f\in {\cal H}^\infty (d).\eqno(3.3)$$
\end{lem}
The proof is the same as that of (2.1) , and we omit it.

   Before continuing, let us observe the following: in view of (2.7), there is a best constant $C_\sigma$ such that, if $f(s)=\sum_1^\infty a_n n^{-s} \in {\cal H}^\infty$ and $\sigma={\Re s}>1/2$, we have :

$$\sum_1^\infty \vert a_n\vert n^{-\sigma}\leq C_\sigma \Vert f\Vert _\infty,\eqno(3.4)$$
 and $C_\sigma\leq (\zeta(2\sigma))^{1/2}$.
We proved in \cite{BCQ} that $C_\sigma =1$ for $\sigma\geq 2$ (this is the analogue of Bohr's inequality in Theorem 1.5), and we might expect that, as in (1.9), as $\sigma\to {1\over 2}$, we have:

$$C_\sigma\sim(\zeta(2\sigma))^{1/2}\sim 1/\sqrt{2\sigma -1},$$ but Theorem 2.5 in Section 2 shows that this is not the case. Now, we have the following Theorem, which shows the optimality of this  Theorem 2.5 and which in particular contains the Bohnenblust -Hille Theorem 2.1  :\\

\begin{theo}: Let $\alpha>0$ be such that
$$\sum_1^\infty  \vert a_n\vert n^{-\sigma_d} (\log n)^\alpha <\infty \textrm{ for each}\  f(s)=\sum_1^\infty a_n n^{-s}\in {\cal H}^\infty(d).$$ Then, $\alpha\leq {d-1\over2}$.
\end{theo}

{\bf Proof 1 (Deterministic):}
We test (3.3) as follows: let $N\geq 2$ be an integer; take the largest integer $s$ such that $p_{s}^d \leq N$, i.e.\ $s=\pi(N^{{1\over d}})$ (recall that $\pi(x)$ is the number of primes $\leq x$), and then take $q=[{\pi(N^{{1\over d}}\over d}]$, so that $r=qd\leq s$. By lemma 3.1, there exists a Dirichlet polynomial $Q(s) =\sum_1^N a_n n^{-s}$ such that $$\Vert Q\Vert_W =q^d\textrm\ {and}\ \Vert Q\Vert_\infty\leq q^{{d+1\over 2}}.$$ Observe that $q\geq C'_d{N^{1\over d}\over \log N}$ by the Prime Number Theorem, and take $f=Q$ in (3.4) to get:
$q^d\leq C_{\beta}{N^{\sigma_d}\over (\log N)^{\beta}} q^{d+1\over 2}$, i.e:
$$(\log N)^\beta \leq C_{\beta}{N^{\sigma_d}\over q^{d-1\over 2}}\leq C_{\beta}C''_{d}{N^{\sigma_d}\over N^{{d-1\over 2d}}}(\log N)^{{d-1\over 2}},$$ or
$(\log N)^\beta\leq C_\beta C''_d(\log N)^{{d-1\over 2}}$, implying $\beta\leq {d-1\over 2}$, and thereby proving Theorem 3.1.

\smallskip
 
{\bf Proof 2 (Probabilistic):}
We use Lemma 2.4 , where we take $y=N^{{1\over d}}, r=\pi(y)$; we denote by $A$ the set of square-free numbers $n$ obtained from $p_1,...,p_r$ and such that $\Omega (n)\leq d$, i.e:
$$n\in A\Leftrightarrow n= p_{i_1}...p_{i_s},\textrm{with}\  1\leq i_1<...<i_s\leq r,\textrm{and}\ s\leq d.$$ 

Let $ a_n =1$  for  $n\in A$, $a_n =0$ for $n\notin A$, and $$P(s)=\sum_{n\in A}a_n n^{-s}\in {\cal H}^\infty (d) \bigcap{\cal P}_N.$$ (indeed, $n\in A\Rightarrow n\leq p_{r}^d\leq y^d =N$).  

 If  $\vert A\vert$ denotes  the cardinality of $A$ and if we test (3.3) on the random polynomial  $P_\omega=\sum_{n\in A}\epsilon_{n}(\omega)a_n n^{-s}$ which is in ${\cal H}^\infty(d)$, we get:\\ 

$$\vert A\vert \leq C_{\beta} N^{\sigma_d} (\log N)^{-\beta} \Vert P_\omega \Vert _\infty.$$

If we now integrate with respect to $\omega$ and use (2.5), noticing that $\Vert P\Vert_2 = \vert A\vert^{1/2}$, we get :
$$\vert A\vert \leq C'_{\beta} N^{\sigma_d}(\log N)^{-\beta}\vert A\vert^{1/2}({N^{1/d}\over \log N})^{1/2}\sqrt{\log\log N},$$ i.e:
$$\vert A\vert^{1/2}\leq C'_{\beta} N^{1/2}(\log N)^{-\beta -1/2}\sqrt{\log\log N}.$$

Now, we  have $\vert A\vert ={r\choose d}\sim r^d/d!$, and we know that
 $r\sim { N^{1/d}\over log(N^{1/d})}$, so $\vert A\vert \sim {N\over(logN)^d} \delta_d$,
and the previous inequality reads:
$${N^{1/2}\over (\log N)^{d/2}}\leq C''_{\beta} N^{1/2}(\log N)^{-\beta -1/2}\sqrt{\log\log N}.$$

This clearly implies:
$-d/2\leq -\beta  -1/2$, i.e $\beta\leq {d-1\over 2}$, proving again (the extra $\log\log$ factor played no role)  Theorem 3.1.
\hfill$\diamondsuit$

We will now allow $q$ and $d$ to vary with $N$, but in such a way that $N^{{1\over d}}\to \infty$, i.e.${\log N\over d}\to \infty$. We will recover the following form of Theorem 2.3:
\begin{theo}: There exists a numerical constant $a>0$ such that:
$$  S(\Lambda_N)\geq a\sqrt N exp(-(1+o(1))\lambda(N))\textrm\  {as}\  N\to \infty.\eqno(3.5)$$
\end{theo}
{\bf Proof (Deterministic):}
Let $d=d_N$ to be adjusted, such that ${\log N\over d}\to \infty$. In Lemma 3.1, we take $$q=[{\pi(N^{{1\over d}}\over d}] \textrm\ {and}\  r=qd,\textrm\ {so\   that } \ p_r^d\leq N,$$
and that the polynomial $Q=Q_d$ of this lemma is in ${\cal P}_N$. We then have by definition:
$$ S(\Lambda_N)\geq {\Vert Q\Vert _W\over \Vert Q\Vert _\infty}\geq {q^d\over q^{{d+1\over 2}}} =q^{{d-1\over 2}}.$$
By the Prime Number Theorem with remainder ($\Pi(x) = {x\over \log x}+O({x\over \log^{2}x}$)), we see that  
$$q\geq {N^{1\over d}\over d \log (N^{1\over d})}(1+O({d\over \log N})={N^{1\over d}\over \log N}(1+O({d\over \log N}),$$ so that 
$$S(\Lambda_N)\geq {N^{d-1\over 2d}\over (\log N)^{d-1\over 2}}(1+O({d^2\over \log N}))$$ $$\geq \sqrt N \exp[-{1\over 2}({\log N\over d}+d\log\log N)](1+O({d^2\over \log N})).$$
We minimize the factor ${\log N\over d} +d\log\log N$ by adjusting 

 $d=[({\log N\over \log\log N})^{1\over 2}]$, so that ${d^2\over\log N}\to 0$. Injecting that value in the preceding inequality gives ( recall that $\lambda (x) = \sqrt {\log x\log\log x}$)
 
$S(\Lambda _N)\geq \sqrt N \exp[ -(1+o(1))\lambda(N)]$, 
and that finishes the proof of Theorem 3.4.\hfill$\diamondsuit$.
 
(We get here   a slightly less good value of $b_0$
 ( $b_0=1$) than the one obtained by R. de la Bret\`eche in \cite{LAB} ($b_0={\sqrt 2\over 2}$), but with a constructive Proof, and using only the PNT, not the delicate estimates connected to Rankin's method).
 It would be interesting to know if we can go farther with Lemma 3.1, but a significant difference should be pointed out: The polynomial $Q$ of lemma 3.1 uses only square-free integers, by construction,  whereas the random constructions ( \cite{Q},\cite{LAB} ) used integers which had only "small" prime factors, but which were not necessarily square-free. It might be the case that the value $b_0=1$ is optimal as long as we use only square-free integers.

\noi{\bf Remark:}
If we do not insist on the best value of $b_0$ in a deterministic proof of Theorem 2.3, we might take for $q$ the largest power of $2$ less than ${\pi(N^{1\over d})\over d}$, and by using only Hadamard matrices we get a minoration of $S(\Lambda _N)$ through polynomials $Q$ whose coefficients are real-valued and indeed take only one of the three values  $+1,-1,0$.

{\bf Acknowledgments:} We thank F.Bayart for indicating us the simple example of the text  in the failure of the Corona Theorem for ${\cal H}^\infty$ on $\C_0$ and for useful discussions. We also thank JP. Kahane for clever remarks on a first version of this paper.

\end{document}